\documentclass[12pt]{article}

\usepackage[english]{babel} 
\usepackage{graphicx}
\usepackage[utf8]{inputenc} 
\usepackage{lmodern}
\usepackage[T1]{fontenc}
\usepackage{url}
\usepackage{xcolor}
\usepackage{sectsty,titling}
\usepackage[top=1.5cm, bottom=2.5cm, outer=2cm, inner=2cm]{geometry} \usepackage{pgf,tikz}
\usetikzlibrary{shapes,positioning,backgrounds}
\usepackage{eso-pic,graphicx}
\usepackage{hyperref}
\usepackage{amsthm}
\usepackage{mathrsfs}
\usetikzlibrary{arrows}
\usepackage{fixltx2e}
\usepackage{microtype}
\usepackage{setspace}
\title{\bf A systematic approach on some relevant theorems that follows from Kolmogorov's axioms}
\author{
Diego J. Raposo \\ 
}
\date{}

\begin{document}

\newtheorem{ax}{Axiom}[section]
\newtheorem{lem}{Lemma}[section]
\newtheorem{thm}{Theorem}[section]
\newtheorem{cor}{Corolary}[section]
\newtheorem{prop}{Proposition}[section]
\newtheorem{defi}{Definition}[section]
\newtheorem{ex}{Example}

\maketitle

\begin{abstract}
\noindent
A selection of the relevant theorems of Probability Theory that comes directly from Kolmogorov's axioms, Set Theory basic results, definitions and rules of inference are listed and proven in a systematic approach, aiming the student who seeks a self-contained account on the matter before moving to more advanced material. 
\end{abstract}

\section{Introduction}
\label{introduction}

Most of the Probability Theory and Statistics books presents the rules of probability as consequences of Andrei Kolmogorov's axioms \cite{ross2010,rozanov1969,magalhaes2006,degroot1989,shiryaev1996,sinai1992,gnedenko1997}. Although they show proofs of the most relevant relations between probabilities of different kinds of events, either directly or through of exercises, I've found no systematic list of relations and proofs. The lacking of some generalizations is also present, so I've made a selection that includes the most common and relevant theorems (and their consequences) that arise, direct or indirectly, from the axioms. The list is neither complete nor fundamentally rigorous, but provides a secure step for the student to base its researches and prove more theorems even before the introduction of random variables. At the end of this paper, in section \ref{probdiagram}, a diagram relating axioms and main results is presented to increase the broad view of the connections among them. 

The discussion is intentionally didactic in order to help the student to follow the reasoning. It demands some contact with proof theory and logic beforehand, but nothing more, alongside some Set Theory equations listed briefly in section \ref{sets}. 

It's perhaps important to emphasize that the axiomatic system proposed by Kolmogorov, somewhat inspired in the frequentist view of statistics \cite{shafer2006}, was not the only one proposed, as described, for instance, by Terenin and Draper \cite{terenin2017}. 

\section{Sets}
\label{sets}

For the proofs to follow, some relations between sets are necessary. I'll present them here, without proof, since they are not the subject of this paper. However, they can be found easily in introductions to mathematical proof \cite{wohlgemuth2011}, for instance, and Probability Theory textbooks:

\begin{itemize}
\item \emph{Empty set}: An empty set, $\emptyset$, have the following properties for any set $A$: $A\cup \emptyset = A$ and $A\cap \emptyset = \emptyset$;
\item \emph{Space set}: An space set, $\Omega$, is the union of all possible sets. In other words, the properties $A \subset \Omega$, $A\cup \Omega = \Omega$ and $A\cap \Omega = A$ are valid for any set $A$;
\item \emph{Complementary set}: A complementary set of $A$, $\overline{A}$, has the following properties: $A \cup \overline{A}= \Omega$ and $A \cap \overline{A}= \emptyset$;
\item \emph{Associative laws for sets}: Given three sets $A$, $B$ and $C$, it can be proven that: $(A\cup B) \cup C = A\cup (B \cup C)$ and $(A\cap B) \cap C = A\cap (B \cap C)$;
\item \emph{Distributive laws for sets}: Given three sets $A$, $B$ and $C$, the following equations are valid: $A\cup(B\cap C) =(A\cup B)\cap (A\cup C)$ and $A\cap(B\cup C) =(A\cap B)\cup (A\cap C)$. 
\end{itemize}

\section{Probabilities}
\label{probabilities}

\begin{defi}
$\Omega$ is the set that represents the sample space, the space of all possible events.
\label{definition1}
\end{defi}

\begin{defi}
Events related to the sample space are all subsets of $\Omega$.
\label{definition2}
\end{defi}

\begin{defi}
Two events $A$ e $B$ are pairwise mutually exclusive (PME) if $A\cap B=\emptyset$, that is, are disjoint sets.
\label{definition3}
\end{defi}

\begin{defi}
A class of subsets of $\Omega$, represented by $\mathcal{F}$, is considered a $\sigma$-algebra if it has the following properties \cite{magalhaes2006,sinai1992}:
\begin{enumerate}
\item $\Omega\in\mathcal{F}$.
\item If $A\in\mathcal{F}$, then $\overline{A}\in\mathcal{F}$.
\item (\emph{Closure with respect to countable unions}). If a countable collection $\{A_1,A_2,...\}=\{A_i\}_{i=1}^{\infty}$ of sets $A_i$ is such that $A_i\in\mathcal{F}$ for all $i$, then $\displaystyle\bigcup_{i=1}^{\infty}A_i\in\mathcal{F}$.
\end{enumerate}
\label{definition4}
\end{defi}

\begin{defi}
A partition of the sample space $\Omega$ is defined according the following property:
\begin{equation}
\displaystyle\bigcup_{i=1}^{n}A_i=\Omega
\label{equation1}
\end{equation}
Being $A_i$ and $A_j$ mutually exclusive (ME) for all combinations of sets.
\label{definition5}
\end{defi}

\begin{defi}
Probability $P$ is a function of the subsets of the sample space correspondent to $\mathcal{F}$: $P=P(\mathcal{F})$. Additionally, it should obey the axioms that follow.
\label{definition6} 
\end{defi}

\begin{center}
\fbox{\parbox{15cm}{
\begin{ax}[Non-negativity]
$P(A_i)\geq 0$, for all event $A_i\subset\mathcal{F}$.
\label{axiom1}
\end{ax}

\begin{ax}[Normalization]
$P(\Omega)=1$.
\label{axiom2}
\end{ax}

\begin{ax}[Countable additivity]
If $A_i$ and $A_j$ are two PME events (disjoint sets) for all $i\neq j$, then $P\displaystyle\left(\bigcup_{i=1}^{\infty}A_i\right)=\displaystyle\sum_{i=1}^{\infty}P(A_i)$ for all $(A_i,A_j)\in\mathcal{F}$.
\label{axiom3}
\end{ax}
} }
\end{center}

\begin{defi}
If $\mathcal{F}$ is a $\sigma$-algebra of set $\Omega$, and $P$ a function of $\mathcal{F}$ with the properties described by axioms \ref{axiom1}, \ref{axiom2} and \ref{axiom3}, then the triple $\{\Omega,\mathcal{F},P\}$ is called \emph{probability space}.
\end{defi}

\begin{thm}
$P(\emptyset)=0$.
\label{theorem1}
\end{thm}

\begin{proof}
According to axiom \ref{axiom3}, we can choose the sets $A_i$ for $i\geq 2$ such that $A_i=\emptyset_i=\emptyset$. Consequently:

\begin{equation}
P\left(\bigcup_{i=1}^{\infty}A_i\right)=P\left[A_1\cup\left(\bigcup_{i=2}^{\infty}A_i\right)\right]=P\left[A_1\cup\left(\bigcup_{i=2}^{\infty}\emptyset_i\right)\right]=P(A_1)+P\left(\bigcup_{i=2}^{\infty}\emptyset_i\right)
\label{equation2}
\end{equation}

Considering that $\displaystyle\bigcup_{i=2}^{\infty}\emptyset_i=\emptyset$, it follows that:

\begin{equation}
P(A_1\cup\emptyset)=P(A_1)+P(\emptyset)
\label{equation3}
\end{equation}

%The second equality being due to axiom \ref{axiom3} again. 
Given the property of the empty set $A_i\cup\emptyset=A_i$ for all $i$, then $P(A_1\cup\emptyset)=P(A_1)$ which can be substituted in Eq. \ref{equation3}:

\begin{equation}
P(A_1)=P(A_1)+P(\emptyset)
\label{equation4}
\end{equation}

Finally proving that:

\begin{equation}
P(\emptyset)=0
\label{equation5}
\end{equation}

\end{proof}

\subsection{Combination of events}
\label{combinationsofevents}

\begin{thm}
If $A_i$ and $A_j$ are PME events, then $P\displaystyle\left(\bigcup_{i=1}^{n}A_i\right)=\displaystyle\sum_{i=1}^{n}P(A_i)$ for $i\neq j$ and $n\geq 1$.
\label{theorem2}
\end{thm}

\begin{proof}
From axiom \ref{axiom3}, consider that from $i=1$ to $i=n$ we have the sets $A_1$, $A_2$, etc, $A_n$, and for $i>n$ we have $A_i=\emptyset_i$. Therefore:

$$P\left(\bigcup_{i=1}^{\infty}A_i\right)=P\left[\left(\bigcup_{i=1}^{n}A_i\right)\cup\left(\bigcup_{i=n+1}^{\infty}A_i\right)\right]=P\left[\left(\bigcup_{i=1}^{n}A_i\right)\cup\left(\bigcup_{i=n+1}^{\infty}\emptyset_i\right)\right]$$

\begin{equation}
P\left[\left(\bigcup_{i=1}^{n}A_i\right)\cup\left(\bigcup_{i=n+1}^{\infty}\emptyset_i\right)\right]=\sum_{i=1}^{n}P(A_i)+\sum_{i=n+1}^{\infty}P(\emptyset_i)
\label{equation6}
\end{equation}

Since $\bigcup_{i=n+1}^{\infty}\emptyset_i=\emptyset$, and $A\cup\emptyset=A$ for all $A$, including $A=\bigcup_{i=1}^{\infty}A_i$, according to theorem \ref{theorem1} $P(\emptyset_i)=P(\emptyset)=0$, thus:

\begin{equation}
P\left[\left(\bigcup_{i=1}^{n}A_i\right)\cup\left(\bigcup_{i=n+1}^{\infty}\emptyset_i\right)\right]=P\left(\bigcup_{i=1}^{n}A_i\right)=\sum_{i=1}^{n}P(A_i)
\label{equation7}
\end{equation}

\end{proof}

\begin{thm}[Normalization condition]
Let the sets $A_1$, $A_2$, ..., $A_n$ be a partition in $\Omega$. It implies that:

\begin{equation}
\displaystyle\sum_{i=1}^nP(A_i)=1
\label{equation8}
\end{equation}

\end{thm}

\begin{proof}
Given the definition \ref{definition5}, $A_i\cap A_j=\emptyset$ for all $i\neq j$. Therefore, following the result of the theorem \ref{theorem2}:

\begin{equation}
\displaystyle\sum_{i=1}^nP(A_i)=P\left(\bigcup_{i=1}^nA_i\right)=P(\Omega)=1
\label{equation9}
\end{equation}

The second equality in Eq. \ref{equation9} follows from definition \ref{definition5}, and the third is a consequence of the axiom \ref{axiom2}.
\end{proof}

\begin{lem}
Any pair among the sets $A_1$, $A_2$, ..., $A_n$ are mutually exclusive (that is, they are PME) if and only if they are mutually exclusive as a whole (ME for any combination of any number of these sets, including all of them). 
\label{lemma1}
\end{lem}

\begin{proof}
If $A_i\cap A_j=\emptyset$ for all $i\neq j$ in the presented sequence of sets, then, given that $A_i\cap\emptyset =\emptyset$:

$$A_1\cap A_2\cap A_3\cap ... \cap A_n=(A_1\cap A_2)\cap A_3\cap ... \cap A_n=\emptyset\cap A_3\cap ... \cap A_n$$
$$A_1\cap A_2\cap A_3\cap ... \cap A_n=(\emptyset\cap A_3)\cap ... \cap A_n=\emptyset\cap ... \cap A_n$$
$$ \ldots $$

\begin{equation}
A_1\cap A_2\cap A_3\cap ... \cap A_n=\emptyset
\label{equation10}
\end{equation}

And since $\bigcap_{i=1}^{n}A_i=\emptyset$, we can deduce that for any $A_i$ and $A_j$ with $i\neq j$:

\begin{equation}
A_1\cap ... \cap A_i ... \cap A_j \cap ... \cap A_n=(A_i\cap A_j)\cap (A_1\cap ... \cap A_n)=\emptyset
\label{equation11}
\end{equation}

Hence, the property $A\cap \emptyset = \emptyset$ for all $A$ can be applied:

$$(A_i\cap A_j)\cap (A_1\cap ... \cap A_n)=\emptyset=\emptyset\cap (A_1\cap ... \cap A_n)$$

\begin{equation}
(A_i\cap A_j)=\emptyset
\label{equation12}
\end{equation}

\end{proof}

\begin{lem}
If $A$ and $B$ are events, then:

\begin{equation}
P(A\cup B)=P(A)+P(B)-P(A\cap B)
\label{equation13}
\end{equation}

\label{lemma2}
\end{lem}

\begin{proof}
By using the relations $A\cup B=A\cup(\overline{A}\cap B)$ and $B=(\overline{A}\cap B)\cup (A\cap B)$, and since $A$ and $\overline{A}\cap B$ are PME, and $\overline{A}\cap B$ and $A\cap B$ too, we can use the theorem \ref{theorem2} for $n=2$ in order to obtain the relations:

\begin{equation}
P(A\cup B)=P(A)+P(\overline{A}\cap B)
\label{equation14}
\end{equation}

\begin{equation}
P(B)=P(\overline{A}\cap B)+P(A\cap B)
\label{equation15}
\end{equation}

Subtracting Eq. \ref{equation15} from \ref{equation14} leads to $P(A\cup B)=P(A)+P(B)-P(A\cap B)$.
\end{proof}

\begin{thm}[Rule of addition of probabilities, or inclusion-exclusion principle, or Poincaré's theorem]

\begin{equation}
P\left(\displaystyle\bigcup_{i=1}^{n}A_i\right)=\sum_{i=1}^{n}P(A_i)-\sum_{i=1}^{n}\sum_{j=i+1}^{n-1}P(A_i\cap A_j)+...
\label{equation16}
\end{equation}

\label{theorem4}
\end{thm}

\begin{proof}
The proof follows from mathematical induction. Eq. \ref{equation16} is refered as proposition $Q(n)$. The case for $n=2$, that is, $Q(2)$, was already proved (lemma \ref{lemma2}), if we change the notation to $A=A_1$ and $B=A_2$. The case $Q(1)$ is trivial, with $P(A)=P(A)$ (or $P(A_1)=P(A_1)$). We can prove the validity of $Q(n+1)$ if the validity of $Q(n)$ is presumed or, equivalently, we might prove $Q(n)$ from $Q(n-1)$. But first we'll prove $Q(3)$ to better understand the structure of $Q(n)$, and of its terms. As already proven:

\begin{equation}
P(A_1\cup A_2)=P(A_1)+P(A_2)-P(A_1\cap A_2)
\label{equation17}
\end{equation}

In order to demonstrate $Q(3)$, we first use $Q(2)$ (lemma \ref{lemma2}) and the associative properties of sets:

\begin{equation}
P(A_1\cup A_2\cup A_3)=P[A_1\cup (A_2\cup A_3)]=P(A_1)+P(A_2\cup A_3)-P[A_1\cap (A_2\cup A_3)]
\label{equation18}
\end{equation}

Then using the distributive property, $A_1\cap (A_2\cup A_3)=(A_1\cap A_2)\cup(A_1\cap A_3)$, and two more applications of $Q(2)$:

$$P(A_2\cup A_3)=P(A_2)+P(A_3)-P(A_2\cap A_3)$$
$$P[A_1\cap (A_2\cup A_3)]=P(A_1\cap A_2)+P(A_1\cap A_3)-P(A_1\cap A_2\cap A_3)$$

\begin{eqnarray}
P(A_1\cup A_2\cup A_3)
& = & P(A_1)+P(A_2)+P(A_3) \nonumber \\
& - & P(A_1\cap A_2)-P(A_1\cap A_3)-P(A_2\cap A_3) \nonumber \\
& + & P(A_1\cap A_2\cap A_3) \nonumber \\
\label{equation19}
\end{eqnarray}

Eq. \ref{equation19} can be written in terms of summations:

\begin{equation}
P\displaystyle\left(\bigcup_{i=1}^{3}A_i\right)=\sum_{i=1}^{3}P(A_i)-\sum_{i=1}^{3-1}\sum_{j=i+1}^{3}P(A_i\cap A_j)+\sum_{i=1}^{3-2}\sum_{j=i+1}^{3-1}\sum_{k=j+1}^{3}P(A_i\cap A_j\cap A_k)
\label{equation20}
\end{equation}

Which for $n$ events can be generalized to:

\begin{eqnarray}
P\displaystyle\left(\bigcup_{i=1}^{n}A_i\right)
& = & \sum_{i=1}^{n}P(A_i)-\sum_{i=1}^{n-1}\sum_{j=i+1}^{n}P(A_i\cap A_j)+\ldots \nonumber \\
& + & (-1)^{L-1}\sum_{i=1}^{n-(L-1)}\sum_{j=i+1}^{n-(L-2)}...\sum_{l=m+1}^{n-(L-L)}P(A_i\cap A_j\cap ...\cap A_m\cap A_l)+\ldots \nonumber \\
& + & (-1)^{n-1}\sum_{i=1}^{n-(n-1)}\sum_{j=i+1}^{n-(n-2)}...\sum_{e=d+1}^{n-(n-n)}P(A_i\cap A_j\cap ...\cap A_d\cap A_e) \nonumber \\
\label{equation21}
\end{eqnarray}

A simpler notation for $Q(n)$ can be:

\begin{eqnarray}
P\displaystyle\left(\bigcup_{i=1}^{n}A_i\right)
& = & \sum_{1\leq i\leq n}P(A_i)-\sum_{1\leq i < j\leq n}P(A_i\cap A_j)+\ldots \nonumber \\
& + & (-1)^{L-1}\sum_{1\leq i < j < ...< m < l\leq n}P(A_i\cap A_j\cap ...\cap A_m\cap A_l) +\ldots \nonumber \\
\label{equation22}
\end{eqnarray}

The equation that corresponds to $Q(n-1)$ can be written as:

\begin{eqnarray}
P\displaystyle\left(\bigcup_{i=2}^{n}A_i\right)
& = & \sum_{2\leq i\leq n}P(A_i)-\sum_{2\leq i < j\leq n}P(A_i\cap A_j)+\ldots \nonumber \\
& + & (-1)^{L-1}\sum_{2\leq i < j < ...< m < l\leq n}P(A_i\cap A_j\cap ...\cap A_m\cap A_l)+\ldots \nonumber \\
\label{equation23}
\end{eqnarray}

And assuming $Q(n-1)$ we should prove $Q(n)$ in order to complete the proof by induction. From the left side of Eq. \ref{equation16}:

\begin{equation}
P\left(\bigcup_{i=1}^{n}A_i\right)=P\left[A_{1}\cup\left(\bigcup_{i=2}^{n}A_i\right)\right]
\label{equation24}
\end{equation}

By applying $Q(2)$:

\begin{eqnarray}
P\left(\bigcup_{i=1}^{n}A_i\right)
& = & P(A_{1})+P\left(\bigcup_{i=2}^{n}A_i\right) - P\left[A_{1}\cap\left(\bigcup_{i=2}^{n}A_i\right)\right] \nonumber \\
& = & P(A_{1})+P\left(\bigcup_{i=2}^{n}A_i\right) - P\left[\bigcup_{i=2}^{n}\left(A_1\cap A_i\right)\right] \nonumber \\
\label{equation25}
\end{eqnarray}

If we apply $Q(n-1)$ in the last two terms on the right of Eq. \ref{equation25}:

\begin{eqnarray}
P\left(\displaystyle\bigcup_{i=1}^{n}A_i\right)
& = & P(A_1)+\displaystyle\sum_{2\leq i\leq n}P(A_i)-\displaystyle\sum_{2\leq i < j\leq n}P(A_i\cap A_j)+\ldots \nonumber \\
& + & (-1)^{L-1}\sum_{2\leq i < j < ...< m < l\leq n}P(A_i\cap A_j\cap ...\cap A_m\cap A_l) +\ldots \nonumber \\
& - & \left\{\displaystyle\sum_{2\leq i \leq n}P(A_1\cap A_i)+\ldots+(-1)^{L-2}\sum_{2\leq i < j < ...< m \leq n}P(A_1\cap A_i\cap ...\cap A_m)+...\right\} \nonumber \\
\label{equation26}
\end{eqnarray}

Performing the substitutions:

\begin{equation}
P(A_1)+\displaystyle\sum_{2\leq i\leq n}P(A_i)=\sum_{1\leq i\leq n}P(A_i)
\label{equation27}
\end{equation}

\begin{equation}
-\displaystyle\sum_{2\leq i \leq n}P(A_1\cap A_i)-\displaystyle\sum_{2\leq i < j\leq n}P(A_i\cap A_j)=-\displaystyle\sum_{1\leq i < j\leq n}P(A_i\cap A_j)
\label{equation28}
\end{equation}

\begin{center}
\ldots
\end{center}

\begin{eqnarray}
\displaystyle -(-1)^{L-2}\sum_{2\leq i < j < ...< m \leq n}P(A_1\cap A_i\cap ...\cap A_m)\nonumber \\ 
+(-1)^{L-1}\sum_{2\leq i < j < ...< m < l\leq n}P(A_i\cap A_j\cap ...\cap A_m\cap A_l)\nonumber \\
=(-1)^{L-1}\sum_{1\leq i < j < ...< m < l\leq n}P(A_i\cap A_j\cap ...\cap A_m\cap A_l) \nonumber \\
\label{equation29}
\end{eqnarray}

We reach $Q(n)$, and the theorem is proved by mathematical induction.

\end{proof}

\begin{lem}
If $A$ and $B$ are mutually exclusive events, then:

\begin{equation}
P(A\cup B)=P(A)+P(B)
\label{equation30}
\end{equation}

\label{lemma3}
\end{lem}

\begin{proof}[Proof 1 (without theorem \ref{theorem2})]
Since the sets are PME, according to definition \ref{definition3} and the theorem \ref{theorem1}: $P(A\cap B)=P(\emptyset)=0$. Considering the lemma \ref{lemma2} the implication is that $P(A\cup B)=P(A)+P(B)-P(A\cap B)=P(A)+P(B)$. 
\end{proof}

\begin{proof}[Proof 2 (by using theorem \ref{theorem2})]
By using the theorem \ref{theorem2} for $n=2$, and changing the notation of $A_1=A$ and $A_2=B$, PME, $P(A\cup B)=P(A)+P(B)$.
\end{proof}

\begin{lem}[Rule of addition of a finite number of ME events]
Let the events $A_1$, $A_2$, \ldots, $A_n$ be ME. The following relation is valid:

\begin{equation}
P\left(\displaystyle\bigcup_{i=1}^{n}A_i\right)=\sum_{i=1}^{n}P(A_i)
\label{equation31}
\end{equation}

\label{lemma4}
\end{lem}

\begin{proof}
Groups of events that are ME as a whole are PME (lemma \ref{lemma4}). Then we can apply the theorem \ref{theorem2} in order to finish the proof.
\end{proof}

\begin{lem}

\begin{equation}
P(\overline{A})=1-P(A)
\label{equation32}
\end{equation}
\label{lemma5}
\end{lem}

\begin{proof}[Proof 1 (based on lemma \ref{lemma3} and two axioms)]
Given that $A$ and $\overline{A}$ are PME, then $P(A\cup\overline{A})=P(A)+P(\overline{A})$, according to axiom \ref{axiom3}. Considering the property of complementary events, $A\cup\overline{A}=\Omega$, we can use the axiom \ref{axiom2}: $P(A\cup\overline{A})=P(\Omega)=1=P(A)+P(\overline{A})$, therefore $P(\overline{A})=1-P(A)$.
\end{proof}

\begin{proof}[Proof 2 (based on \ref{lemma2} and one axiom)] 

Assume $B=\overline{A}$, then $P(A\cup\overline{A})=P(A)+P(\overline{A})-P(A\cap\overline{A})$. Since $P(A\cup\overline{A})=P(\Omega)$, according to the axiom \ref{axiom2}, we have $P(A\cup\overline{A})=1$, and given the fact that $A$ and $\overline{A}$ are PME (another property of the complementary sets), $P(A\cap\overline{A})=0$ (definition \ref{definition3}). Therefore, $1=P(A)+P(\overline{A})$, leading to Eq. \ref{equation32}.
\end{proof}

\begin{lem}
If $A\subset B$, then $P(A)\leq P(B)$.
\label{lemma6}
\end{lem}

\begin{proof}
If $A\subset B$, then $B=A\cup(\overline{A}\cap B)$. Being $A$ and $\overline{A}\cap B$ PME, then lemma \ref{lemma3} applies, and $P(B)=P(A)+P(\overline{A}\cap B)$. Following the axiom \ref{axiom1}, $P(\overline{A}\cap B)\geq 0$, then $P(B)\geq P(A)$ and the proof is complete.
\end{proof}

\begin{lem}

\begin{equation}
P(A)\leq 1
\label{equation33}
\end{equation}

\label{lemma7}
\end{lem}

\begin{proof}
By the definition of the space set, $A\subset\Omega$, and from lemma \ref{lemma6} $P(A)\leq P(\Omega)$. The direct application of the axiom \ref{axiom2}, leads to the result: $P(A)\leq P(\Omega)=1$.
\end{proof}

\subsection{Dependency among events}
\label{dependencyevents}

\begin{defi}
The conditional probability of the event $A$ given the event $B$, $P(A|B)$, is defined for $P(B)>0$ as:

\begin{equation}
P(A|B)=\displaystyle\frac{P(A\cap B)}{P(B)}
\label{equation34}
\end{equation}

\label{definition8}
\end{defi}

\begin{lem}
Given $A$ and $B$, and $P(B) > 0$, it's true that $0\leq P(A|B)\leq 1$.
\label{lemma8}
\end{lem}

\begin{proof}
According the definition of $P(A|B)$, that depends on $P(A\cap B)$ and $P(B)>0$, once $0\leq P(C)\leq 1$ for any $C$ (axiom \ref{axiom1} and lemma \ref{lemma7}), then both $0\leq P(A|B)\leq 1$ and $0\leq P(A\cap B) \leq 1$ are true. The first condition is the proof of the lemma, and the second can be used to determine the limits of $P(A\cap B)$, leading to $0\leq P(A\cap B)\leq P(B)$. When $P(A\cap B)=P(B)$, $P(A|B)=1$. And if $P(A\cap B)=0$, that is, the events $A$ and $B$ are PME, then $P(A|B)=0$. 
\end{proof}

\begin{prop}
If $A$ and $B$ are PME events, then $P(A|B)=0$.
\label{proposition1}
\end{prop}

\begin{proof}
If $A$ and $B$ are PME, $A\cap B=\emptyset$ (definition \ref{definition3}), and since $P(\emptyset)=0$ (theorem \ref{theorem1}), then $P(A|B)=0$. Based on the definition of $P(A|B)$ (equation \ref{equation34}), we have $P(A|B)=0$.
\end{proof}

\begin{prop}
If the event $B$ implies the event $A$, that is, $B\subset A$, then $P(A|B)=1$.
\label{proposition2}
\end{prop}

\begin{proof}
If $B\subset A$, $P(A\cap B)=P(B)$. By definition, $P(A|B)=P(A\cap B)/P(B)$, hence $P(A|B)=P(B)/P(B)=1$. 
\end{proof}

\begin{prop}[Rule of addition for conditional probabilities]
If $A_1$, $A_2$, ..., $A_n$ are ME with union $A=\displaystyle\bigcup_{i=1}^nA_i$, then:

\begin{equation}
P(A|B)=\displaystyle\sum_{i=1}^nP(A_i|B)
\label{equation35}
\end{equation}

\label{proposition3}
\end{prop}

\begin{proof}
According the definition of $P(A|B)$ (definition \ref{definition8}) and its relation with $P(B)>0$ and $P(A\cap B)$:

\begin{equation}
P(A|B)=\displaystyle\frac{P(A\cap B)}{P(B)}=\frac{1}{P(B)}P\displaystyle\left(B\cap\bigcup_{i=1}^{n}A_i\right)=\frac{1}{P(B)}P\displaystyle\left[\bigcup_{i=1}^{n}(A_i\cap B)\right]
\label{equation36}
\end{equation}

Since $A_i$ and $A_j$ are PME for $i\neq j$, then $A_i\cap B$ and $A_j\cap B$ also are PME:

\begin{equation}
A_i\cap A_j=\emptyset \Rightarrow (A_i\cap A_j)\cap B=\emptyset\cap B =\emptyset\Rightarrow (A_i\cap B)\cap (A_j\cap B)=\emptyset
\label{equation37}
\end{equation}

Therefore we can apply lemma \ref{lemma4}:

\begin{equation}
P(A|B)=\displaystyle\frac{1}{P(B)}P\displaystyle\left[\bigcup_{i=1}^{n}(A_i\cap B)\right]=\displaystyle\frac{1}{P(B)}\sum_{i=1}^{n}P(A_i\cap B)=\sum_{i=1}^{n}\frac{P(A_i\cap B)}{P(B)}=\displaystyle\sum_{i=1}^{n}P(A_i|B)
\label{equation38}
\end{equation}

\end{proof}

\begin{defi}
The event $A$ is said to be independent of the event $B$, or statistically independent (SI), if and only if $P(A|B)=P(A)$.
\label{definition9}
\end{defi}

\begin{thm}[Rule of the product of probabilities]
Let the events $A_1$, $A_2$, ..., $A_n$ with $P\left(\displaystyle\bigcup_{i=1}^{n}A_i\right)\geq 0$ and $P\left(\displaystyle\bigcup_{i=1}^{k}A_i\right)> 0$ for $1\leq k < n$. It can be shown that:

\begin{equation}
P\left(\displaystyle\bigcap_{i=1}^{n}A_i\right)=P(A_1)\times P(A_2|A_1)\times \ldots \times P(A_n|A_1\cap ...\cap A_{n-1})
\label{equation39}
\end{equation}

\label{theorem5}
\end{thm}

\begin{proof}
Using the definition of $P(A|B)$ repeatedly for each factor in the product:

$$P\displaystyle\left(\bigcap_{i=1}^{n}A_i\right)=\frac{P\displaystyle\left(\bigcap_{i=1}^{n}A_i\right)}{P\displaystyle\left(\bigcap_{i=1}^{n-1}A_i\right)}\times\frac{P\displaystyle\left(\bigcap_{i=1}^{n-1}A_i\right)}{P\displaystyle\left(\bigcap_{i=1}^{n-2}A_i\right)}\times ...\times \frac{P(A_1\cap A_2)}{P(A_1)}P(A_1)$$

\begin{equation}
P\displaystyle\left(\bigcap_{i=1}^{n}A_i\right)=P\left[A_n|\displaystyle\left(\bigcap_{i=1}^{n-1}A_i\right)\right]\times P\left[A_{n-1}|\displaystyle\left(\bigcap_{i=1}^{n-2}A_i\right)\right]\times ... \times P(A_2|A_1)P(A_1)
\label{equation40}
\end{equation}

\end{proof}

\begin{lem}
If $A$ and $B$ are SI, and both $P(A)$ and $P(B)$ are not zero, then $P(A\cap B)=P(A)P(B)$
\label{lemma9}
\end{lem}

\begin{proof}
The definition for $P(A|B)=P(A\cap B)/P(B)$ also applies backwards, $P(B|A)=P(A\cap B)/P(A)$. If $A$ and $B$ are SI, given the definition \ref{definition9} $P(A|B)=P(A)$ and $P(B|A)=P(B)$. Thus, in both cases one can show that $P(A\cap B)=P(A)P(B)$, since $P(A\cap B)=P(A|B)P(B)=P(A)P(B)$ and $P(A\cap B)=P(B|A)P(A)=P(B)P(A)$.
\end{proof}

\begin{defi}
The events $A_1$, $A_2$, ..., $A_n$ are defined as mutually independents (MI) if $P\displaystyle\left(\bigcap_{i=1}^{n}A_i\right)=\prod_{i=1}^{n}P(A_i)$ for all combinations of sets between $1$ e $n$.
\label{definition10}
\end{defi}

\begin{lem}
If the events $A_1$, $A_2$, ..., $A_n$ are MI, then any pair $A_i$ e $A_j$, for $i\neq j$, are SI. 
\label{lemma10}
\end{lem}

\begin{proof}
According the definition \ref{definition10} $A_i$ and $A_j$ must be SI, for the independence is valid for the combination of any number of sets, including pairs. Hence $P(A_i\cap A_j)=P(A_i)P(A_j)$ for any $i\neq j$ if the events $A_1$, $A_2$, ..., $A_n$ are MI. Notice that the opposite is not necessarily true: by assuming $P(A_i\cap A_j)=P(A_i)P(A_j)$ for any pair $(i,j)$ one do not prove $P(A_i\cap A_j\cap A_k)=P(A_i)P(A_j)P(A_k)$ for all $(i,j,k)$, and higher order groups. 
\end{proof}

\begin{lem}
If the sets $A$ and $B$ are SI, then they are not PME, and vice versa.
\label{lemma11}
\end{lem}

\begin{proof}
If $A$ and $B$ are SI, then $P(A\cap B)=P(A)P(B)$ (lemma \ref{lemma9}). PME events are such that $A\cap B=\emptyset$, and $P(A\cap B)=P(\emptyset)$. According theorem $\ref{theorem1}$, $P(\emptyset)=0$, thus PME events ($P(A\cap B)=0$) cannot be SI ($P(A\cap B)=P(A)P(B)$). Notice that the events $A$ and $B$ are such that $P(A)>0$ and $P(B)>0$, respectively (lemma \ref{lemma9}).
\end{proof}

\begin{lem}
If the set $\{C_1,...,C_n\}$ is a partition of the sample space, $\Omega$, then for any event $A$:

\begin{equation}
P(A)=\displaystyle\sum_{i=1}^{n}P(A|C_i)P(C_i)
\label{equation41}
\end{equation}

\label{lemma12}
\end{lem}

\begin{proof}
The definition of $P(A|C_i)$ (Eq. \ref{equation34}) implies $P(A|C_i)=P(A\cap C_i)/P(C_i)$. Making the summation from $i=1$ to $n$ on both sides of the equation $P(A|C_i)P(C_i)=P(A\cap C_i)$:

\begin{equation}
\displaystyle\sum_{i=1}^nP(A|C_i)P(C_i)=\sum_{i=1}^nP(A\cap C_i)
\label{equation42}
\end{equation}

Since $C_i$ and $C_j$ are disjoint for any $i\neq j$, then is also true that $(A\cap C_i)\cap (A\cap C_j)=\emptyset$ (see proof of proposition \ref{proposition3}). Therefore, by using the theorem \ref{theorem2}:

\begin{eqnarray}
\displaystyle\sum_{i=1}^nP(A|C_i)P(C_i)
& = & P\left[\bigcup_{i=1}^n (A\cap C_i)\right]=P\left[A\cap \left(\bigcup_{i=1}^n C_i\right)\right]\nonumber \\
& = & P(A\cap \Omega) = P(A|\Omega)P(\Omega)=P(A) \nonumber \\
\label{equation43}
\end{eqnarray}

Where we have used the definition of a partition in $\Omega$ (definition \ref{definition5}), the axiom \ref{axiom2}, the definition of conditional probability and the implicit identity $P(A|\Omega)$ for any $A$ ($A\subset \Omega$, so any element of the set $A$ is also an element of the sample space).

\end{proof}

\begin{thm}[Bayes' theorem]

Given the event $A$ is such that $P(A)\geq 0$, and the set $\{C_1, ...,C_n\}$, which defines a partition in $\Omega$, with $P(C_i)\geq 0$ for every $i$, it's possible to prove that:

\begin{equation}
P(C_i|A)=\displaystyle\frac{P(A|C_i)P(C_i)}{\displaystyle\sum_{i=1}^{n}P(A|C_i)P(C_i)}
\label{equation44}
\end{equation}

\label{theorem6}
\end{thm}

\begin{proof}
Due to the fact that $P(C_i|A)=P(A\cap C_i)/P(A)$ and $P(A|C_i)=P(A\cap C_i)/P(C_i)$ (definition \ref{definition8}), then:

\begin{equation}
P(C_i|A)=\frac{P(A|C_i)P(C_i)}{P(A)}
\label{equation45}
\end{equation}

We can substitute the result of lemma \ref{lemma12} for $P(A)$ in Eq. \ref{equation45} to prove the theorem. 
\end{proof}

\subsection{Probabillity theorems diagram}
\label{probdiagram}

\begin{figure}[!htbp]
\centering
\includegraphics[width=\textwidth]{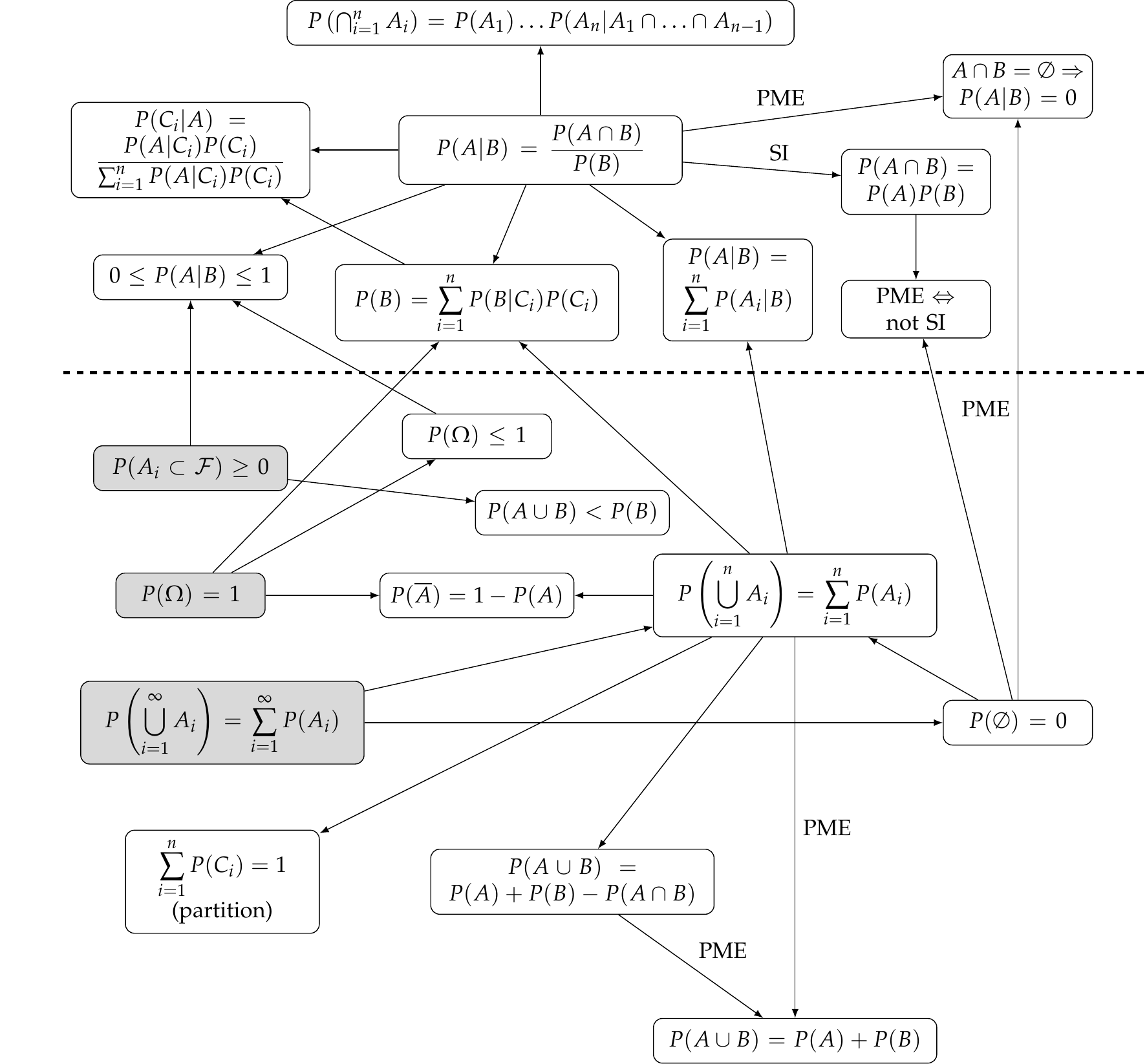}
\caption{Main results from Kolmogorov's axioms (in gray boxes). Those above the dashed line correspond to dependent events, and those below the line are related to combinations of events.}
\label{diagram}
\end{figure}

In order to represent the main content of this paper, namely, the most relevant equations presented so far, and the relations among them, the diagram in Fig. \ref{diagram} was prepared. It omits some results and assumptions, focusing on the fundamental relations between axioms and results. A line divides the figure in two parts: below are the results proved in section \ref{combinationsofevents}, and above it the results from section \ref{dependencyevents} are listed and interrelated.

\section{Conclusion}
\label{conclusion}

I hope the set of proofs and the choice of theorems/lemmas/propositions and definitions, as much as the order they are presented, help students of mathematics, statistics, engineering, chemistry and physics to se a broad picture of the Kolmogorov's axiomatic system, even if in a simplified and incomplete form. 

\section{Acknowledgements}
\label{acknowledgements}

I must thank CNPq and CAPES for the PhD scholarship which allowed me to investigate the theme.


\begin{thebibliography}{99}

\bibitem{degroot1989} DeGroot, M. H. Probability and Statistics. Addison-Wesley Publishing Company, 1989.

\bibitem{gnedenko1997} Gnedenko, B. V. Theory of Probability, 6th ed. CRC Press, 1997.

\bibitem{magalhaes2006} Magalhães, M. N. Probabilidade e Variáveis Aleatórias, 2nd ed. Edusp, 2006.

\bibitem{ross2010} Ross, S. A First Course in Probability, 8th ed. Prentice Hall, 2010.

\bibitem{rozanov1969} Rozanov, Y. A. Probability Theory: A Concise Course. Dover Publications, 1969.

\bibitem{shafer2006} Shafer, G., and Vovk, V. The sources of kolmogorov’s grundbegriffe. Statistical Science 21, 1 (2006), 70–98.

\bibitem{shiryaev1996} Shiryaev, A. Probability, 2nd ed. Springer, 1996.

\bibitem{sinai1992} Sinai, Y. Probability Theory: An Introdutory Course. Springer, 1992.

\bibitem{terenin2017} Terenin, A., and Draper, D. Cox’s theorem and the jaynesian interpretation of probability. \url{https://arxiv.org/abs/1507.06597v2} (2017), 1–18.

\bibitem{wohlgemuth2011} Wohlgemuth, A. Introduction to Proof in Abstract Mathematics. Dover Publications, 2011.

\end{thebibliography}
\end{document}